\documentstyle[11pt]{article}
\def\Bbb R{{\rm \bf R}}
\def\proclaim#1{\vskip2mm{\bf #1}\em}
\def\endproclaim{\em \vskip2mm}
\def\tag#1{\eqno(#1)}
\def\gathered{\begin{array}{c}}
\def\endgathered{\end{array}}
\def\text{\mbox}

\begin{document}

\title {Remarks on Lin-Nakamura-Wang's paper}
\author{Masaru IKEHATA\footnote{
Laboratory of Mathematics,
Graduate School of Engineering,
Hiroshima University, Higashihiroshima 739-8527, JAPAN}}
\maketitle

\begin{abstract}
Theorem 1.2 in their paper arXiv:1904.00999v1 [math.AP] 30 Mar 2019
``Reconstruction of unknown cavity by single measurement'' is not valid.


\noindent
AMS: 35R30

\end{abstract}


\section{A counter example}

In \cite{LNW} they state\footnote{Please refer to their paper \cite{LNW} for the symbols used in this note without explanation.}
if $\overline D\not\subset\overline G$, then $I(G)=\infty$.
However, in this note we give a simple example that $\displaystyle\overline D\not\subset\overline G$, however $I(G)=0$.

Let $\Omega=\{x\in\Bbb R^2\vert\,\vert x\vert<R\}$ with $R>1$ and $D=\{x\in\Bbb R^2\,\vert\,\vert x\vert<1\}$.
Let $u$ solve
$$
\left\{\displaystyle
\begin{array}{ll}
\displaystyle
\Delta u=0 & \text{in $\Omega\setminus\overline D$,}
\\
\\
\displaystyle
\frac{\partial u}{\partial\nu}=0 & \text{on $\partial D$,}\\
\\
\displaystyle
u(R\cos\theta,R\sin\theta)=\left(R+\frac{1}{R}\right)\cos\theta, & \text{$\theta\in\,[0,\,2\pi[$}.
\end{array}
\right.
\tag {1.0}
$$
Note that the solution has the explict form
$$\displaystyle
u(r\cos\theta,r\sin\theta)=\left(r+\frac{1}{r}\right)\cos\theta.
$$
The key point of this note is the following trivial fact:
$u$ has an extension to the domain $\tilde{\Omega}=\{x\in\Bbb R^2\vert\,0<\vert x\vert<R\}=\Omega\setminus\{0\}$ as a solution of the Laplace equation.

Let $0<\delta<1$ and choose $G=\{x\in\Bbb R^2\,\vert\,\vert x\vert<1-\delta\}$.
We have $\overline G\subset D$ and thus $\overline D\not\subset\overline G$.

Given $\epsilon>0$ let $g\in H^{1/2}(\partial\Omega)$ be an arbitrary function
such that the solution $z_g$ of 
$$\left\{
\begin{array}{ll}
\displaystyle
\Delta z_g=0 & \text{in $\Omega$,}\\
\\
\displaystyle
z_g=g & \text{on $\partial\Omega$}
\end{array}
\right.
$$
satisfies
$$\displaystyle
\Vert z_g\Vert_{H^1(G)}<\epsilon.
\tag {1.1}
$$

By Lemma 2.1 in \cite{LNW} we have
$$\displaystyle
\int_{\partial\Omega}\partial_{\nu}w\cdot g\,ds
=-\int_{\partial D} u\cdot\partial_{\nu}z_g\,ds,
\tag {1.2}
$$
where $w=u-v$ and $v$ solves
$$\left\{
\begin{array}{ll}
\displaystyle
\Delta v=0 & \text{in $\Omega$,}\\
\\
\displaystyle
v=u & \text{on $\partial\Omega$.}
\end{array}
\right.
$$
Let $\tilde{u}$ denote the {\it harmonic extension} of $u$ into $\tilde{\Omega}$, that is
$$\displaystyle
\tilde{u}(r\cos\theta,r\sin\theta)=\left(r+\frac{1}{r}\right)\cos\theta.
$$
Let $C=\{x\in\Bbb R^2\,\vert\,\vert x\vert=1-\delta'\}$ with $\delta<\delta'<1$.
We have $C\subset G$.

Write
$$\begin{array}{l}
\displaystyle
\,\,\,\,\,\,
-\int_{\partial D} u\cdot\partial_{\nu}z_g\,ds
\\
\\
\displaystyle
=\int_{\partial D}\left(\partial_{\nu}u\cdot z_g-u\cdot\partial_{\nu}z_g\right)\,ds
\\
\\
\displaystyle
=\int_{\partial D}\left(\partial_{\nu}\tilde{u}\cdot z_g-\tilde{u}\cdot\partial_{\nu}z_g\right)\,ds.
\end{array}
$$
Since $\tilde{u}$ and $z_g$ are harmonic in $1-\delta'<\vert x\vert<1$, one has the expression
$$\displaystyle
\int_{\partial D}\left(\partial_{\nu}\tilde{u}\cdot z_g-\tilde{u}\cdot\partial_{\nu}z_g\right)\,ds
=\int_{C}\left(\partial_{\nu}\tilde{u}\cdot z_g-\tilde{u}\cdot\partial_{\nu}z_g\right)\,ds.
$$
Thus (1.2) becomes
$$\displaystyle
\int_{\partial\Omega}\partial_{\nu}w\cdot g\,ds
=\int_{C}\left(\partial_{\nu}\tilde{u}\cdot z_g-\tilde{u}\cdot\partial_{\nu}z_g\right)\,ds.
$$
It is easy to see that this right-hand side has the bound $O(\Vert z_g\Vert_{H^1(G)})$.
Thus the condition (1.1) yields
$$\displaystyle
\left\vert\int_{\partial\Omega}\partial_{\nu}w\cdot g\,ds
\right\vert
\le C\epsilon,
$$
where $C$ is independent of $g$.  Hence $I_{\epsilon}(G)\le C\epsilon$ and $I(G)=\lim_{\epsilon\downarrow 0}I_{\epsilon}(G)=0$.

\section{Looking at the example in Section 1 a little more}

Let $u$ be the solution of (1.0) and $\tilde{u}$ its harmonic extension to $\tilde{\Omega}$.
In this section $G$ denotes an arbitrary open subset of $\Omega$ such that
$\overline G\subset\Omega$ and $\Omega\setminus\overline G$ is connected.
In this section we prove

\proclaim{\noindent Proposition 2.1.}

(a)  If $(0,0)\in G$, then $I(G)=0$. 

(b)  If $(0,0)\not\in\overline G$, then, for all $\epsilon$ $I_{\epsilon}(G)=\infty$.

\endproclaim

{\it\noindent Proof.}  First we prove (a).  In this case one can find a cirecle $S$ centered at $(0,0)$ such that
$S\subset G$.  At this time, the following equation is obtained as in the previous section:
$$\displaystyle
\int_{\partial\Omega}\partial_{\nu}w\cdot g\,ds
=\int_{S}\left(\partial_{\nu}\tilde{u}\cdot z_g-\tilde{u}\cdot\partial_{\nu}z_g\right)\,ds.
$$
Note that $z_g$ is the same as before.  Thus this together with (1.2) yield $I_{\epsilon}(G)\le C\epsilon$ with a positive constant $C$
independent of $g$.  And hence $I(G)=\lim_{\epsilon\downarrow 0}I_{\epsilon}(G)=0$.

Next we prove (b).  For this we claim the identity:
$$\displaystyle
\int_{\partial\Omega}\partial_{\nu}w\cdot g\,ds
=-2\pi\nabla z_g(0,0)\cdot\mbox{\boldmath $e$}_1,
\tag {2.1}
$$
where $\mbox{\boldmath $e$}_1=(1,0)^T$.

First of all admit equation (2.1) and move on.  
Consider the case $(0,0)\not\in\overline G$.
One can find an open disc $B$ centered at $(0,0)$ and radius $t_0$ such that $\overline B\subset\Omega\setminus\overline G$.
Let $B_t=\{x\in\Bbb R^2\,\vert\,\vert x\vert<t\}$ with $0<t<t_0$.
Since the function
$$\displaystyle
E_t(x)=\log \vert x-t\mbox{\boldmath $e$}_1\vert
$$
is harmonic in a neighbourhood of $\overline{G}\cup \overline{B_{t/2}}$, the Runge approximation property yields: there exists a sequence $\{g_j\}$ such that
$$\displaystyle
\lim_{j\rightarrow\infty}\Vert z_{g_j}-E_t\Vert_{H^1(G\cup B_{t/2})}=0.
\tag {2.2}
$$
Then an interior regulerity estimate yields $z_{g_j}$ together with its all derivatives converges to $E_t$ 
and the corresponding derivatives compact uniformly in $B_{t/2}$.   Thus (2.1) yields
$$\displaystyle
\lim_{j\rightarrow\infty}\int_{\partial\Omega}\partial_{\nu}w\cdot g_j\,ds
=\frac{2\pi}{t}.
\tag {2.3}
$$
Note also that we have 
$$\displaystyle
\lim_{j\rightarrow\infty}\Vert z_{g_j}\Vert_{H^1(G)}=\Vert E_t\Vert_{H^1(G)}.
$$

Given $\epsilon>0$ define
$$\displaystyle
\tilde{g}_j=\frac{\epsilon}{2\Vert E_t\Vert_{H^1(G)}}g_j.
$$
Since the map $g\mapsto z_g$ is linear, we have
$$\displaystyle
\Vert z_{\tilde{g}_j}\Vert_{H^1(G)}=\frac{\epsilon}{2\Vert E_t\Vert_{H^1(G)}}\Vert z_{g_j}\Vert_{H^1(G)}<\epsilon
$$
for all $j>>1$. 

And (2.3) gives
$$\displaystyle
\lim_{j\rightarrow\infty}\int_{\partial\Omega}\partial_{\nu}w\cdot \tilde{g}_j\,ds
=\frac{2\pi}{t}\cdot\frac{\epsilon}{2\Vert E_t\Vert_{H^1(G)}}
\tag {2.4}
$$
Since $\overline B\cap\overline G=\emptyset$, Lebesgue's dominated convergence theorem gives $\lim_{t\downarrow 0}\Vert E_t\Vert_{H^1(G)}=\Vert E_0\Vert_{H^2(G)}<\infty$.  Thus the right-hand side on (2.4) blows up as $t\downarrow 0$.  This yields $I_{\epsilon}(G)=\infty$.

\noindent
$\Box$

{\bf\noindent Remarks.}

(i)  The case $(0,0)\in\partial G$ seems delicate (at the present time).

(ii) This type of sequence satisfying (2.2) has been used in the {\it probe method} \cite{IProbe} which aims at reconstructing unknown
discontinuities such as cavities, inclusions and cracks.  However, the probe method employs the Dirichlet-to-Neumann map,
i.e., infinitely many pairs of the Cauchy data of the governing equation.  Instead in the proof of (b) a single pair of
Cauchy data is {\it fixed} and sequences $z_{g_j}$ produced by {\it infinitely many} $g_j$ are used as test functions.

(iii)  The choices of $\{g_j\}$ in two cases (a) and (b) are different.  Since we do not know the position of $\{(0,0)\}$ in advance,
we have the question: what is the {\it good choice} of $\{g_j\}$ {\it common} to two cases.  
This is also a problem about the no response test.

\subsection{Proof of (2.1)}

Same as before, we have, for all circles $S_{\eta}$ centered at $(0,0)$ with radius $\eta\in\,]0,\,1[$
$$\displaystyle
\int_{\partial\Omega}\partial_{\nu}w\cdot g\,ds
=\int_{S_{\eta}}\left(\partial_{\nu}\tilde{u}\cdot z_g-\tilde{u}\cdot\partial_{\nu}z_g\right)\,ds.
$$
We compute the limt of this right-hand side as $\eta\downarrow 0$.

First we have
$$\begin{array}{l}
\displaystyle
\,\,\,\,\,\,
\int_{S_{\eta}}\partial_{\nu}\tilde{u}\cdot z_g\,ds
\\
\\
\displaystyle
=\left(1-\frac{1}{\eta^2}\right)\eta
\int_0^{2\pi}
\cos\theta\cdot z_g(\eta\cos\theta,\eta\sin\theta)d\theta\\
\\
\displaystyle
=-\left(1-\frac{1}{\eta^2}\right)\eta\int_0^{2\pi}
\sin\theta\cdot \frac{d}{d\theta}\left\{z_g(\eta\cos\theta,\eta\sin\theta)\right\}d\theta\\
\\
\displaystyle
=-\left(1-\frac{1}{\eta^2}\right)\eta^2
\int_0^{2\pi}
\sin\theta\cdot \nabla z_g(\eta\cos\theta,\eta\sin\theta)
\cdot(-\sin\theta,\cos\theta)^T\,d\theta
\\
\\
\displaystyle
\rightarrow
\int_0^{2\pi}
\sin\theta\cdot \nabla z_g(0,0)
\cdot(-\sin\theta,\cos\theta)^T\,d\theta
\\
\\
\displaystyle
=-\pi\nabla z_g(0,0)\cdot\mbox{\boldmath $e$}_1.
\end{array}
$$
Second we have
$$\begin{array}{l}
\displaystyle
\,\,\,\,\,\,
\int_{S_{\eta}}\tilde{u}\cdot\partial_{\nu}z_g\,ds
\\
\\
\displaystyle
=(\eta^2+1)\int_0^{2\pi}\cos\theta\cdot\nabla z_g(\eta\cos\theta,\eta\sin\theta)\cdot
(\cos\theta,\sin\theta)^T\,d\theta
\\
\\
\displaystyle
\rightarrow
\pi\nabla z_g(0,0)\cdot\mbox{\boldmath $e$}_1.
\end{array}
$$
This completes the proof.

\section{One can not apply Fatou' s lemma}

The key point of their argument on page 5 is the definiteness of the signature of
$\partial_{\nu_x}F_{\mbox{\boldmath $a$}}(x,y)$ for $x\in N_{y_0}\cap\partial D$ and $y\rightarrow y_0$ along the axis of the cylinder
$N_{y_0}$.  Here we give an example of $D$ that does not ensure this property.

Let $D$ be a bounded domain and in $x_3<0$.
We assume that $y_0=(0,0,0)\in\partial D$ and $N_{y_0}\cap\partial D$ is {\it flat}
and included in the plane $x_3=0$.
Thus $\nu_x=\nu_{y_0}=\mbox{\boldmath $e$}_3$.  

Let $E(x)=\frac{1}{\vert x\vert}$.
We have
$$\displaystyle
\partial_3E(x)=-\frac{x_3}{\vert x\vert^3},
$$
and
$$\displaystyle
\partial_3^2E(x)
=\frac{1}{\vert x\vert^5}
(3x_3^2-\vert x\vert^2).
$$

Since $\mbox{\boldmath $a$}=\nu_{y_0}=\mbox{\boldmath $e$}_3$, 
we have, for all $x\in N_{y_0}\cap\partial D$ and $y=(0,0,y_3)$ with $0<y_3<<1$
$$\displaystyle
\partial_{\nu_x}F_{\mbox{\boldmath $a$}}(x,y)=-\partial_3^2E(x-y)
$$
and thus
$$\displaystyle
\partial_{\nu_x}F_{\mbox{\boldmath $a$}}(x,y)=
-\frac{1}{\vert x-y\vert^5}
(2y_3^2-x_1^2-x_2^2).
$$
Therefore we have

(i) if $x_1^2+x_2^2<2y_3^2$, then $\partial_{\nu_x}F_{\mbox{\boldmath $a$}}(x,y)<0$;

(ii) if $x_1^2+x_2^2>2y_3^2$, then $\partial_{\nu_x}F_{\mbox{\boldmath $a$}}(x,y)>0$.

Thus as $y_3\downarrow 0$ the sign of the function $\partial_{\nu_x}F_{\mbox{\boldmath $a$}}(x,y)$ of
$x\in N_{y_0}\cap\partial D$ can not have a definite sign.

This implies, one can not apply Fatou's lemma as done (3.4) in this simplest case.

\section{Another reason of invalidness of (3.5) on page 5: A heuristic explanation}

Even general case one can not obtain (3.5).
Its heuristic explanation is the following.

Since $\mbox{\boldmath $a$}=\nu_{y_0}$, if $x\in N_{y_0}\cap \partial D$
we expect
$$\displaystyle
\partial_{\nu_x}F_{\mbox{\boldmath $a$}}(x,y)
\sim
-\partial_{\nu_{x_0}}^2E(x-y).
$$
However, $E$ satisfies the Laplace equation we have
$$\displaystyle
\partial_{\nu_{x_0}}^2E(x-y)
=-(\partial_{x_1}^2+\partial_{x_2}^2)E(x-y),
$$
where $x_1$ and $x_2$ are {\it tangential directions} at $y_0$.  Thus we can expect
$$\displaystyle
\partial_{\nu_x}F_{\mbox{\boldmath $a$}}(x,y)
\sim
(\partial_{x_1}^2+\partial_{x_2}^2)E(x-y).
$$
Then the integral
$$\displaystyle
\int_{N_{y_0}\cap\partial D}u(x)\cdot\partial_{\nu_x}F_{\mbox{\boldmath $a$}}(x,y)ds(x)
$$
may become
$$\displaystyle
\sim
\int_{N_{y_0}\cap\partial D}u(x)\cdot
(\partial_{x_1}^2+\partial_{x_2}^2)E(x-y)ds(x).
$$
Then applying integration by parts to this right-hand, one can reduce the singularity of integrand twice and gets 
an integral and additional terms which are bounded as $y\rightarrow y_0$.

\section{Some comments on references}

In \cite{I1} (1999!) using a single set of the Cauchy data, we have already given the reconstruction formula of the convex hull of
unknown polygonal cavity $D$ and done its numerical testing in \cite{IO}. 
The method developed in this paper is called the {\it enclosure method} and based on the asymptotic behaviour
of the integral with respect to a large parameter $\tau$
$$\displaystyle
\int_{\partial\Omega}\partial_{\nu}w\,g\,ds,
$$
where $g=e^{\tau x\cdot(\omega+i\omega^{\perp})}$ with two unit vectors $\omega$ and $\omega^{\perp}$ 
perpendicular each other.  Note that in this case $z_g(x)=e^{\tau x\cdot(\omega+i\omega^{\perp})}$.

Besides, in the case when $\Omega$ is an ellipse, even though the homogeneous background is {\it unknown},
the enclosure method works and yields
a reconstruction formula of the convex hull of the union of the polygonal cavity and the focal points of $\Omega$ by using 
a single flux corresponding to a {\it band-limited} surface potential \cite{I3}.

These informations are missed in \cite{LNW}.

\section{Extendability}

The point is the extendability of the potential $u$ from $\Omega\setminus\overline D$ across $\partial D$ into $D$,
for example, if $\partial D$ is a real analytic surface, then by applying the Cauchy-Kovalevskaya theorem one has
such an extension locally. 
In this case, we can prove that, by doing the procedure above locally around $y_0\in\partial D\setminus\overline G$ on page 5 in \cite{LNW},
(3.5) in \cite{LNW} is not valid.  The enclosure method in \cite{I1} catches a corner where one can not have an extention 
of the potential (due to Friedman-Isakov's extension argument \cite{FI} under the condition
$\text{diam}\,D<\text{dist}\,(D,\partial\Omega)$).

So at least we have to find an argument that employs explicitly the impossibility of applying the Cauchy-Kovalevskaya theorem
on $\partial D$.

\section{Conclusion}

The problem is not simple and still unsolved! 
I guess the complete version of the no response test with a single measurement
tells us the limt of the extension of the soultion (continuation as a solution of the governing equation).
Proposition 2.1 is an evidence of this belief.

$$\quad$$

\centerline{{\bf Acknowledgments}}

The author was partially supported by Grant-in-Aid for
Scientific Research (C)(No. 17K05331) and (B)(No. 18H01126) of Japan  Society for
the Promotion of Science.

$$\quad$$

\vskip1cm
\noindent
e-mail address

ikehata@hiroshima-u.ac.jp

\end{document}